\documentclass[12pt,reqno]{amsart}
\usepackage{color}
\usepackage{latexsym}
\usepackage{amsmath,amssymb,dsfont,amsfonts}
\usepackage{mathrsfs}
\usepackage{verbatim}
\usepackage{graphicx}
\usepackage{graphics}
\usepackage{geometry}
\usepackage{booktabs}
\usepackage[shortlabels]{enumitem}
\usepackage{xcolor}
\newtheorem{remark}{Remark}[section]
\newtheorem{proposition}[remark]{Proposition}
\newtheorem{theorem}[remark]{Theorem}
\newtheorem{definition}[remark]{Definition}
\newtheorem{corollary}[remark]{Corollary}
\newtheorem{lemma}[remark]{Lemma}
\newtheorem{assumption}[remark]{Assumption}

\def\R{\mathbb R}

\def\E{\mathbb E}
\def\P{\mathbb P}
\def\Q{\mathbb Q}

\def\shd{{\mathcal D}}
\def\shf{{\mathcal F}}

\def\shl{{\mathcal L}}

\newenvironment{prooff}{{\bf \textit{Proof}}}{\hfill $\Box$ \\}
\newenvironment{preuve}{{\bf \textit{Proof.}}}{\hfill $\Box$ \\}

\numberwithin{equation}{section}

\allowdisplaybreaks

\title{
On path-dependent SDE\MakeLowercase{s} involving distributional
drifts}
\author{Alberto OHASHI$^1$}

\author{Francesco RUSSO$^2$}

\author{Alan TEIXEIRA$^1$}

\address{$1$ Departamento de Matem\'atica, Universidade de Bras\'ilia, 70910-900, Bras\'ilia, Brazil.} \email{amfohashi@gmail.com}

\address{$2$ ENSTA Paris, Institut Polytechnique de Paris,
 Unit\'e de Math\'ematiques appliqu\'ees, 828, boulevard des Mar\'echaux, F-91120 Palaiseau, France}
 \email{francesco.russo@ensta-paris.fr}
\address{$3$ ENSTA Paris, Institut Polytechnique de Paris,
 Unit\'e de Math\'ematiques appliqu\'ees, 828, boulevard des Mar\'echaux, F-91120 Palaiseau, France}
\email{alan.teixeira@ensta-paris.fr}

\date{April 25th 2021}

\begin{document}

\begin{abstract}
 In this paper, we study (strong and weak) 
existence and uniqueness of a class of
 non-Markovian SDEs whose drift contains the derivative in the sense of distributions
of a continuous function.
\end{abstract}

\maketitle



\section{Introduction}

\label{SIntro}

This paper discusses in detail a framework of
one-dimensional stochastic differential equations (henceforth abbreviated by SDEs) with distributional drift
with possible path-dependency.
Even though we could have worked in the multidimensional case,
we have preferred to explore in a systematic way the real line case.

The main objective of this paper is to analyze the solution 
(existence and uniqueness) of the  martingale problem associated with SDEs of the type
\begin{equation}\label{X}
dX_t = \sigma(X_t)dW_t + b'(X_t)dt + \Gamma(t, X^t)dt,\ \ X_0 \stackrel{d}{=}
\delta_{x_0},
\end{equation}
where $b, \sigma: \mathbb{R} \rightarrow \mathbb{R}$ are continuous
functions, $\sigma > 0, x_0 \in \R$
and $W$ is a standard Brownian motion.
The assumptions on $b$, which will be formulated later, imply that $b'$ is a Schwartz distribution.
Concerning the path-dependent component of the drift, we consider a locally bounded
functional
\begin{equation} \label{EGamma}
\Gamma: \Lambda \rightarrow \mathbb{R},
\end{equation}
where 
\begin{equation*}
\Lambda:= \{(s, \eta) \in [0, T] \times C([0, T]); \eta = \eta^{s}\}
\end{equation*}
and
\begin{equation*}
\eta^{s}(t)=\left\{\begin{array}{lr}
	\eta(t), & \text{if} \ t \leq s\\
	\eta(s), & \text{if} \ t > s.
	\end{array}\right.
\end{equation*}
By convention, we extend $\Gamma$ from $\Lambda$ to $[0,T] \times C([0,T])$
by setting (in non-anticipating way) 
$$ \Gamma(t,\eta):= \Gamma(t,\eta^t), t \in [0,T], \eta \in C([0,T]).$$
Path-dependent SDEs were investigated under several aspects.
Under standard Lipschitz regularity conditions on the coefficients, it is known, (see e.g Theorem 11.2 \cite[chapter V]{rogers})  
that strong existence and uniqueness holds.
In case the path-dependence takes the form of delayed stochastic equations, 
one-sided Lipschitz condition ensures strong existence and uniqueness,
see e.g \cite{scheutzow,mohammed1984stochastic}.
Beyond Lipschitz regularity on the coefficients of the SDE, \cite{frikha}
shows  uniqueness in law under structural conditions on an underlying approximating Markov 
process, where local-time and running maximum dependence are
 considered. Weak existence of infinite-dimensional SDEs with additive noise on 
the configuration space with path-dependent drift functionals with sublinear growth are studied by \cite{dereudre}. In all these works, the drift is a non-anticipative functional.
Beyond Brownian motion based driving noises,  \cite{castrequini}
 establishes existence of solutions for path-dependent Young differential
 equation.

 The Markovian case ($\Gamma = 0$) with
distributional drift  has been intensively studied over the years.
Diffusions in the generalized sense were
first considered in the case when the solution is still a semimartingale,
 beginning with, at least in our
knowledge  \cite{portenko}. Later on, many
authors considered special cases of SDEs with generalized coefficients. It is
difficult to quote them all; in the case when the drift $b'$ is a measure
and the solutions are semimartingales,  we refer
the reader to  \cite{blei, esd, trutnau2}.
 We also recall that
 \cite{ew} considered even special cases of non-semimartingales
solving stochastic differential equations with generalized drift; 
those cases include examples coming from Bessel processes, for which only existence is proved.

In the non-semimartingale  case, time-independent SDEs in dimension one
of the type
\begin{equation}  \label{SDEone}
\mathrm dX_t = \sigma(X_t) \mathrm dW_t +  b'(X_t)\mathrm dt,\quad t\in[0,T],
\end{equation}
where $\sigma$ is a strictly positive continuous
function and $b'$ is the derivative of
a real continuous function
was solved and analyzed carefully in
 \cite{frw2} and \cite{frw1}, which  treated well-posedness
 of the martingale   problem, It\^o's formula under weak conditions,
 semimartingale
 characterization and Lyons-Zheng
decomposition.
The only supplementary assumption was the
 existence of the function
\begin{equation}\label{Sigmaintr}
\Sigma (x) :=  2\int_{0}^{x}\frac{b'}{\sigma^{2}}(y)dy, x \in \mathbb{R},
\end{equation}
considered as a suitable limit via regularizations.
Those authors considered weak solutions. The SDE \eqref{SDEone} was also investigated
 by \cite{basschen1}, where the authors  provided a well-stated framework when $\sigma$
 and $b$ are $\gamma$-H\"older
 continuous,  $ \gamma > \frac{1}{2}$. In \cite{russo_trutnau07},
 the authors have
 also shown that in some cases strong solutions
 exist and pathwise uniqueness holds.
 More recently, in the time-dependent framework (but still one-dimensional),
 a significant contribution
was done by \cite{Delarue2017ROUGHPA}.
As far as the multidimensional case is concerned,
some important steps were done in \cite{flandoli_et.al14} and
more recently in \cite{cannizzaro2018},
when the diffusion matrix is the identity and $b$ is a time-dependent
drift in some suitable negative Sobolev space.
We also refer to   \cite{basschen2}, where the authors
have focused on \eqref{X}
in the case of a time independent drift $b$ which
is a measure of Kato class.
To complement the list, we study existence and uniqueness of a class of
path-dependent SDEs whose drift contains the derivative in the sense of
distributions
of a continuous function. To our best knowledge, this is the first
paper which approaches a class of non-Markovian SDEs with distributional
drifts.

In this work, equation \eqref{X} will be interpreted as a martingale problem
 with respect to some operator  $\shl f := L f + \Gamma f',$ 
see \eqref{DOperatorL},  
where $L$ is the Markovian generator 
\begin{equation} \label{DefL}
Lf = \frac{\sigma^{2}}{2}f'' + b'f',
\end{equation}
where we stress that $b'$ is the derivative of some continuous
function $b$, 
 interpreted in the sense of distributions.  
If we  denote $\Sigma$ as in \eqref{Sigmaintr}, then
the operator $L$ can be written as
 \begin{equation}
 \label{DefLbis}
Lf = (e^\Sigma f')'\frac{e^{-\Sigma}\sigma^2}{2},
\end{equation}
see \cite{frw1}. We define a notion of martingale problem related to $\shl$ (see Definition \ref{D31})
and a notion of strong martingale problem related to $\shd_\shl$ and a given Brownian
motion $W$, see Definition \ref{DSolution}; that definition
 has to be compared with the notion of strong existence and pathwise uniqueness 
 of an SDE.
 In the Markovian case, the notion of {\it strong martingale problem}
 was introduced in
\cite{russo_trutnau07} and it represents the corresponding notion
to strong solution of SDEs in the framework of martingale problems.

As anticipated, we will concentrate on the case when $b$ is continuous, the case
of special discontinuous functions is investigated in
\cite{ORT1_Bessel}. When $\Gamma = 0$, this case was
completely analyzed in \cite{frw1} and \cite{frw2}. Concerning the former case, under the existence of the function (\ref{Sigmaintr}) and some boundedness 
or linear growth condition (see \eqref{H}), Theorem \ref{P41} 
presents existence for  the martingale problem related to \eqref{X}.
Proposition \ref{uni} states   uniqueness under more 
restrictive conditions.
Under suitable Lipschitz regularity conditions 
on a functional $\tilde{\Gamma}$, which is related to $\Gamma$ via \eqref{GammaTilde}, 
Corollary \ref{pauniMart}
establishes well-posedness for the
strong martingale problem associated to \eqref{X}.

That type of process appears in some fields for instance, much work has been done on
Markovian processes in a random environment. In this particular case, $b = B$ is a two-sided 
real-valued Brownian motion which is independent from $W$ and (\ref{X})
might be interpreted 
as 
the non-Markovian version 
$$ dX_t = -\frac{1}{2}\dot{B}(X_t)dt + \Gamma(t,X^t)dt + dW_t,$$
of the so-called Brox diffusion
which is indeed obtained setting $\Gamma = 0$,
see e.g \cite{brox,yorhushi} and other references therein in the classical Markovian
context.

The paper is organized as follows. After this Introduction
and after having fixed some preliminairies, in Section \ref{SNonMark}
we define the suitable concept of martingale problems for 
SDE with distributional drift with path-dependent perturbations.
In Section \ref{SDerCont} we investigate the case when $b'$ 
is the derivative of a continuous function.

\section{Notations and Preliminaries }

\label{SNotations}
\subsection{General notations}
$ $

\setcounter{equation}{0}

Let $I$ be an interval of $\mathbb{R}$.
$C^{k}(I)$ is the space of real functions  defined on $I$ having
 continuous derivatives till order $k$.
  Such space is endowed with the uniform convergence topology on compact sets
for the functions and all derivatives.
Generally, $ I = \mathbb{R}$ or $[0, T]$ for some fixed
positive real $T$.
The space of continuous functions on $I$ will be denoted by $C(I)$.
Often, if there is no ambiguity $C^{k}(\mathbb{R})$ will be simply indicated by
$C^{k}.$
Given an a.e. bounded real function $f$, $\vert f \vert_\infty$
will denote the essential supremum.

We recall some notions from \cite{frw1}.
For us, all filtrations $\mathfrak{F}$ fulfill the usual conditions.
When no filtration is specified, we mean the canonical filtration of the
underlying process.
Otherwise the canonical filtration associated with a process $X$
is denoted by $\mathfrak{F}^X$.
An $\mathfrak{F}$-Dirichlet process $X$ is the sum of an $\mathfrak{F}$-local martingale $M^X$ with an $\mathfrak{F}$-adapted zero quadratic variation process $A^X$.
We will fix by convention that $A_0^X = 0$ so that the decomposition is unique.
A sequence $(X^n)$ of continuous processes indexed by $[0,T]$ is said to converge u.c.p.
to some process $X$ whenever $\displaystyle\sup_{t \in [0,T]} \vert X^n_t - X_t \vert$ converges to zero in probability.

\vfill \eject
\begin{remark} \label{RDir}
	$ $
\begin{enumerate}
\item An $\mathfrak{F}$-continuous semimartingale $Y$ is
  always an $\mathfrak{F}$-Dirichlet process. The $A^Y$ process
coincides with the continuous bounded variation component.
Moreover the quadratic variation $[Y]$ is the usual quadratic
variation for semimartingales.
\item Any $\mathfrak{F}$-Dirichlet process is a finite quadratic variation
process
and its quadratic variation gives $[X] = [M^X]$.
\item If $f \in C^1(\R)$ and $X =M^{X} + A^{X}$ is an $\mathfrak{F}$-Dirichlet process,
then  $Y:= f(X)$ is again an $\mathfrak{F}$-Dirichlet process and
$[Y] = \int_0^\cdot f'(X)^2 d[M^{X}]$.
\end{enumerate}
\end{remark}

\section{Non-Markovian SDE: the function case. }\label{SNonMark}

\subsection{General considerations.}

$ $

Similarly as for the case of Markovian SDEs, it is possible to formulate
the notions of strong existence, pathwise uniqueness, existence and uniqueness
in law for path-dependent SDEs of the type \eqref{X}, see
e.g. Section \ref{Appendix}.

Let us suppose for the moment that
$\sigma, b': \R \rightarrow \R$ are Borel functions.
We will consider solutions $X$ of
\begin{equation} \label{E1.3}
\left \{
\begin{array}{cll}
dX_t  &=& \sigma(X_t)  dW_t + b'(X_t)dt + \Gamma(t,X^t) dt   \\
X_0  & = & \xi,
\end{array}
\right.
\end{equation}
for some initial condition $\xi$.
The previous equation will be denoted by $E(\sigma,b', \Gamma; \nu)$
(where $\nu$ is the law of $\xi$),
or simply with $E(\sigma,b', \Gamma)$ if we omit the initial condition.
For simplicity in this paper, $\xi$  will always be considered
as deterministic.
$\Gamma$ is defined in \eqref{EGamma}.

 \begin{definition} \label{D112Bis}
Let $\nu$ be the Dirac  probability measure on $\R$ such that
$\nu = \delta_{x_0}, x_0 \in \R$.
A stochastic process $X$ 
 is called {\bf (weak) solution} 
of  $E(\sigma,b', \Gamma;\nu)$
with respect to a probability $\P$ if there is a Brownian motion $W$
on some filtered probabilty space,
such that $X$ solves \eqref{E1.3}  and $X_0 = x_0.$
We also say that the couple $(X,\P)$ solves $E(\sigma,b', \Gamma)$
with initial condition distributed according to $\nu$.
\end{definition}

Suppose $\Gamma \equiv 0$.
A very well-known result in \cite{sv}, Corollary 8.1.7, concerns the
equivalence between martingale problems and solution in law.
Suppose for a moment  that $b'$ is a continuous function.
According to \cite[chapter 5]{karatzasShreve}, it is well-known
that a process $X$ and probability $\P$ solve the classical martingale problem, if and only if, $X$ is a (weak) solution of \eqref{X}.
The proof of the result mentioned above can be easily adapted to the path-dependent case, i.e.
 when
$\Gamma \neq 0$. This provides
 the statement below.
\begin{proposition} \label{P31}
A couple $(X,\P)$ is a solution of
$E(\sigma, b',\Gamma)$,
if and only if, under $\P$, 
\begin{equation}\label{Ef32}
f(X_t) - f(X_0) - \int_0^t Lf(X_s) ds - \int_0^t f'(X_s) \Gamma(s,X^s)ds
\end{equation}
is a local martingale, where
$ Lf = \frac{1}{2} \sigma^2 f'' + b'f'$,
for every $f \in C^2$.
\end{proposition}

\subsection{Comments about the distributional case}

\label{C1ter}

$ $

When $b'$ is a distribution,
it is not obvious to introduce the notion of SDE, except in
the case  when $L$ is close to the divergence form, i.e. when
$Lf = (\sigma^2 f')' + \beta f'$ and $\beta$ is a Radon measure,
see e.g. Proposition 3.1 of \cite{frw1}.
For this reason, we replace the notion
of weak solution with the notion of  martingale problem.
Suppose for a moment that $L$ is a second order PDE operator
with possible generalized coefficients.
In general, as it is shown in \cite{frw1},
$C^2$ is not included in the domain of operator $L$ and,
similarly to
 \cite{frw1}, we will replace $C^2$ with some domain $\mathcal{D}_L$.
Suppose that $L: \shd_L \subset C^1(\R) \rightarrow C(\R)$.
 

\begin{definition} \label{D31}
\begin{enumerate}	
\item We say that a continuous stochastic process $X$ solves 
 (with respect to a probability $\P$ on some measurable space 
$(\Omega,\shf)$)
the
martingale problem related to
\begin{equation} \label{DOperatorL}
\shl f := L f + \Gamma f',
\end {equation}
 with initial condition $\nu = \delta_{x_0}, x_0 \in \R$, 
with respect to a domain $\shd_L$
	if
	\begin{equation}\label{lmp}
	M^f_t := f(X_t) - f(x_0) - \int_{0}^{t}\mathit{L}f(X_s)ds - \int_{0}^{t}f'(X_s)\Gamma(s, X^s)ds,
	\end{equation}	
	is a $\P$-local martingale for all $f \in \mathcal{D}_{\mathit{L}}$.

We will also say that the couple $(X,\P)$ is a solution of (or $(X,\P)$ solves) the
martingale problem with respect to $\shd_L$.
\item If a solution exists, we say that the martingale problem above
{\it admits existence}.
\item We  say that the martingale problem above
{\it admits uniqueness} if any two solutions $(X^i,\P^i), i  = 1,2$
(on some measurable space $(\Omega,\shf)$)
have the same law.

\end{enumerate}

\end{definition}
In the sequel, when the measurable space $(\Omega,\shf)$ 
is self-explanatory it will be often omitted.
As already observed in Proposition \ref{P31},
the notion of martingale problem is (since the works of Stroock and Varadhan \cite{sv}) a concept related to solutions of SDEs in law.
In the case when $b'$ and $\sigma$ are continuous functions (see \cite{sv}),
$\shd_L$ corresponds to $C^2(\R)$.

Below we introduce the analogous notion of strong existence
and pathwise uniqueness for our martingale problem.
\vfill \eject

\begin{definition} \label{DSolution}
$ $ 

\begin{enumerate}
\item 
Let $(\Omega, \mathcal{F}, \P)$ be a probability space and let
  $ \mathfrak{F} = (\shf_t)$ be the canonical filtration associated with a 
fixed Brownian motion $W$.
Let $x_0 \in \R$ be a constant.
We say that a continuous $\mathfrak{F}$-adapted real-valued process $X$
such that 
$X_0 = x_0$ is a {\bf solution to 
the strong martingale problem}
(related to \eqref{DOperatorL})
 with respect to $\shd_L$ and $W$ (with related filtered probability space),
if
\begin{equation}\label{lmpBis}
	f(X_t) - f(x_0) - \int_{0}^{t}\mathit{L}f(X_s)ds - \int_{0}^{t}f'(X_s)\Gamma(s, X^s)ds = \int_{0}^{t}f'(X_s) \sigma(X_s) dW_s,
	\end{equation}	
	is an $ \mathfrak{F}$-local martingale
for all $f \in \shd_L$.
\item  
We say that the martingale problem 
 related to \eqref{DOperatorL}
 with respect to $\shd_L$ 
admits \textbf{strong existence}
 if for every $x_0 \in \R$, 
 given a filtered probability space
 $(\Omega, \mathcal{F}, \P, \mathfrak{F}),$
where  $\mathfrak{F} = (\shf_t)$ is the canonical filtration associated with 
 a Brownian motion $W$, 
there is a  process $X$ solving the strong martingale problem
(related to \eqref{DOperatorL})
 with respect to $\shd_L$ 
 and $W$
 with $X_0 = x_0$.
\item 
We say that the martingale problem  
(related to \eqref{DOperatorL})
with respect to $\shd_L$,
admits \textbf{pathwise uniqueness} 
if given
$  (\Omega, \mathcal{F}, \P),$
a Brownian motion $W$ on it,
two solutions $X^i, i = 1, 2$, of the strong martingale with respect to
 $\shd_L$ and $W$
  and $\P[X_0^1 = X_0^2] = 1,$ 
then $X^1$ and $X^2$ are indistinguishable.
\end{enumerate}

\end{definition}

\section{The case when the drift is a derivative of
a continuous function}

\label{SDerCont}

 

Here we extend the
Markovian framework of \cite{frw1}
 to the non-Markovian case.

\subsection{The Markovian case}
$ $

\label{S41}

Let $\sigma$ and $b$ be functions in  $C(\R)$ with $\sigma > 0$.
In \cite{frw1}, the authors define, by mollification methods, the function
\begin{equation} \label{ESigma}
 \Sigma(x) = 2 \lim_{n\rightarrow \infty} \int_0^x  \frac{b_n'}
{\sigma_n^2}(y) dy, \forall x \in \R,
\end{equation}
where the limit is intended to be in $C(\R)$, i.e. uniformly on each compact.
For concrete examples, we refer the reader to \cite{frw1},
for instance if either $\sigma^2$ or $b$ are of locally bounded variation.
The function $Lf$ is defined according to \eqref{DefLbis}.

Similarly as in \cite{frw1}, $\shd_L$ will be the linear space of $f \in C^1(\R)$
for which there exists $l \in C(\R)$ such that $Lf = l$
(still by mollifications). Thereby we use
the terminology $Lf = l$ in the $C^1$-generalized sense.
As in  \cite{frw1}, we define the $L$-{\it harmonic} function 
$h: \R \rightarrow \R$ as 
\begin{equation}
\label{eqh}
h(0) = 0, \quad h' = e^{-\Sigma},
\end{equation}
in particular $Lh= 0$, as we will see in Corollary \ref{C22}.
In this case, $L: \shd_L \subset C^1 \rightarrow \R$
can be written as in \eqref{DefLbis}.
Indeed, Proposition \ref{L22} below
is a direct consequence of  Lemma 2.9  and Lemma 2.6 in \cite{frw1}.
From now on $h$ will be the function defined in \eqref{eqh}.

\begin{proposition} \label{L22}
	$ $
\begin{enumerate}
\item Let $f \in C^1$.
	$f \in \mathcal{D}_L$ if, and only if, there exists $\phi \in C^1$ such that $f' = \exp(-\Sigma)\phi$.
\item If 	$f \in \mathcal{D}_L,$ then $l = Lf$ is given by
\eqref{DefLbis}; in particular  we have
\begin{equation} \label{ESquare}
Lf = \phi'  \exp(-\Sigma) \frac{\sigma^2}{2},
\end{equation}
where $\phi$ is the function given in item (1) above.
\end{enumerate}
\end{proposition}

\begin{corollary} \label{C22}
	$ $
\begin{enumerate}
\item If $f \in \shd_L,$ then $f^2 \in \shd_L$ and
$L f^2 = \sigma^2 f'^2 + 2 f Lf.$
\item 
$$ Lh = 0, L h^2 = \sigma^2 h'^2.$$
\end{enumerate}
\end{corollary}
\begin{prooff} \ (of Corollary \ref{C22}). 
\begin{enumerate}
\item $f^2 \in \shd_L$ because $(f^2)' = 2ff' = (2f\phi)\exp(-\Sigma)$ taking into account Proposition \ref{L22} (1) and the fact that $\phi_2 := 2f\phi \in C^1$. By \eqref{ESquare}, 
\begin{eqnarray*}
Lf^2 &=& \phi'_{2}\exp(-\Sigma)\dfrac{\sigma^2}{2} = (f\phi)'\exp(-\Sigma)\sigma^2  \\
&=& f'\sigma^2 \exp(-\Sigma)\phi + f\phi' \exp(-\Sigma)\sigma^2 = f'^2 \sigma^2 + 2fLf.
\end{eqnarray*}
\item It follows by Proposition \ref{L22}, setting $\phi = 1$ and item (1).
\end{enumerate}
\end{prooff}
We now formulate a standing assumption.
\begin{assumption} \label{A1}
	$ $
\begin{itemize}
\item $\Sigma$ is well-defined, see \eqref{ESigma}.
\item We suppose the {\it non-explosion condition}
\begin{equation} \label{Noexpl}
 \int_{-\infty}^0 e^{-\Sigma}(x) dx =  \int_0^{+\infty}  e^{-\Sigma}(x) dx
= \infty.
\end{equation}

\end{itemize}
\end{assumption}

\begin{remark} \label{RNE}
	$ $
\begin{enumerate}
\item Under Assumption \ref{A1},
 the $L$-{\it harmonic} function $h:\R \rightarrow \R$ 
defined in \eqref{eqh} is a
$C^1$-diffeomorphism. In particular $h$ is surjective.
\item  
  It is easy to verify that Assumption \ref{A1} implies
 the non-explosion condition (3.16)   in Proposition
3.13 in \cite{frw1} is fulfilled. 
\end{enumerate}
\end{remark}

\begin{remark}\label{A}
  When $\sigma$ and $b'$ are continuous functions, then $\shd_{L} = C^2$.
  Indeed, in this manner, $\Sigma \in C^1$ and then $f' = \exp(-\Sigma)\phi \in C^1$. In particular, $L f$ corresponds to its classical definition.
\end{remark}

In relation to the {\it harmonic} function $h$ defined in \eqref{eqh},
Proposition 2.13 in \cite{frw1} states the following.

\begin{proposition} \label{P22}
	$f \in \mathcal{D}_L$ if and only if $f\circ h^{-1} \in C^2$.
\end{proposition}
In fact previous result can be generalized to the case when
$h$ is replaced by a bijection $g$
of class $C^1$.



\begin{proposition} \label{Ph}
We suppose Assumption \ref{A1}.
Let $g \in \shd_L$ be a diffeomorphism of class $C^1$ such that $g' > 0$ and let $f$ be a $C^1$ function. Then, $f \in \mathcal{D}_L$ if and only if $f \circ g^{-1}$ belongs to
 $\mathcal{D}_{L^Y}$, where
 $L^Y v := \frac{1}{2} (\sigma^g_0)^2 v'' + ((L g) \circ g^{-1})  v', v \in \mathcal{D}_{L^Y}$,
 and
\begin{equation} \label{Esigma0}
\sigma^g_{0} := (\sigma g') \circ g^{-1}.
\end{equation}
Moreover
$$ L f \circ g^{-1} = L^Y (f \circ g^{-1}).$$

\end{proposition}
\begin{remark} \label{AB}
  By Remark \ref{A}, we get
  $\shd_{L^Y} = C^2$, since $L^Y$ has continuous coefficients.
   \end{remark}

\begin{prooff} \ (of Proposition \ref{Ph}).
By Proposition \ref{L22} there exists $\phi^{g} \in C^1$ such that
\begin{equation} \label{EPhih}
  g' = \exp(-\Sigma)\phi^{g}.
\end{equation}  
Concerning the direct implication, if $f \in \mathcal{D}_{L}$, first we prove that $f\circ g^{-1} \in \mathcal{D}_{L^Y}$. Again by
Proposition \ref{L22} there exists $\phi^f \in C^1$ such that
 $ f' = \exp(-\Sigma)\phi^f.$
So, $(f\circ g^{-1})' = \dfrac{f'}{g'} \circ g^{-1} =
\dfrac{\phi^f}{\phi^{g}} \circ g^{-1} \in C^1$ because $g^{-1} \in C^1$ and $\phi^{g} > 0$, so $f \circ g^{-1} \in C^2$. Note that, by Remark \ref{AB}, $\mathcal{D}_{L^Y} = C^2$ so $f \circ g^{-1} \in \mathcal{D}_{L^Y}$. Moreover, by Proposition
 \ref{L22} (2), \begin{equation*}
(\phi^{f})' = \dfrac{2Lf}{\sigma^2}\exp(\Sigma), \ \ \ (\phi^{g})' = \dfrac{2Lg}{\sigma^2}\exp(\Sigma).
\end{equation*}

A direct computation gives
\begin{equation*}
(f\circ g^{-1})'' = \left(\dfrac{\phi^f}{\phi^g}\circ g^{-1}\right)' = \left[\dfrac{2Lf}{g'^2 \sigma^2} - \dfrac{2Lg f'}{\sigma^2 g'^3}\right]\circ g^{-1}.
\end{equation*}
Consequently,
\begin{equation}\label{ig}
\dfrac{(\sigma_{0}^g)^2}{2}(f \circ g^{-1})'' = \left[Lf - \dfrac{Lg f'}{g'}\right]\circ g^{-1}.
\end{equation}

By \eqref{ig}, $Lf \circ g^{-1} = \dfrac{(\sigma^g_{0})^2}{2}(f \circ g^{-1})'' + (Lg)\circ g^{-1}(f \circ g^{-1})' = L^{Y}(f \circ g^{-1})$.

Let us discuss the converse implication.
Suppose that $f \circ g^{-1}$ belongs to
$\mathcal{D}_{L^Y} = C^2$. 
By Proposition \ref{L22} we need to show that $f' \exp(\Sigma) \in C^1$
which is equivalent to show that
$(f' \exp(\Sigma)) \circ g^{-1}$ belongs to $C^1$.
If $\phi^g \in C^1$ is such that $g' = \exp(-\Sigma)\phi^{g}$ (see \eqref{EPhih})
we have
$$ (f' \exp(\Sigma)) \circ g^{-1} = (f' \frac{\phi^g}{g'}) \circ g^{-1}
= (f \circ g^{-1})' (\phi^g \circ g^{-1}),$$
which obviously belongs to $C^1$. Therefore $f \in \shd_L$.
\end{prooff}
\begin{remark}   \label{Remh}
If
$f \in \shd_L$,
setting $\varphi = f \circ h^{-1}$ in Proposition
\ref{Ph} (with $g=h$),
gives
$$ L (\varphi  \circ h) \circ h^{-1}   = L^Y (\varphi) = 
\frac{1}{2}\sigma_0^2 \varphi'',$$
where
\begin{equation}\label{Esigma0true}
  \sigma_0 = (\sigma h') \circ h^{-1}.
  \end{equation}
\end{remark}
In    \cite{frw1}, the
 authors also show that the existence and uniqueness of the solution of
 the martingale problem
 are conditioned to a non-explosion feature.
An easy consequence of Proposition 3.13 in \cite{frw1}
gives the following.
\begin{proposition} \label{P21}
Let $\nu = \delta_{x_0}, x_0 \in \R $.
Suppose  the validity of Assumption \ref{A1}.
Then 
 the martingale problem 
related to $L$ (i.e. with $\Gamma = 0$)
with respect to $\shd_L$ with initial condition $\nu$
admits existence and uniqueness.
\end{proposition}
\begin{remark} \label{P32frw1}
  By Proposition 3.2 of \cite{frw1}, if $\Gamma = 0$ and $(X,\P)$
  is a solution of the above-mentioned martingale problem, then
there exists a $\P$-Brownian motion $W$ such that \eqref{lmp}
equals 
$$ \int_0^t (f'\sigma)(X_s) dW_s, t  \in [0,T].$$

\end{remark}



\subsection{The path-dependent case.}
\label{SPD}

$ $

Let $\sigma$ and $b$ be functions in $C(\R)$ with $\sigma > 0$ and 
$\Gamma$ as defined in \eqref{EGamma}.
Let us suppose again Assumption  \ref{A1}
and let $h$ be the function defined in \eqref{eqh}.
We recall that $\sigma_0$ was defined in \eqref{Esigma0true}.

The first result explains how to reduce our
path-dependent martingale problem to a path-dependent SDE.

\begin{proposition}\label{X-Y}
Let $(\Omega, \shf, \P)$ be a probability space.
Let $X$ be a stochastic process, we denote $Y= h(X)$.
\begin{enumerate}
\item	$(X, \P)$
solves
  the martingale problem 
related to \eqref{DOperatorL}
with respect to $\shd_L$ if and only if the process $Y := h(X)$ 
 is a solution (with respect to $\P$) of
\begin{equation}\label{Y}
Y_t = Y_0 + \int_{0}^{t}\sigma_0(Y_s)dW_s + \int_{0}^{t}h'(h^{-1}(Y_s)
)\Gamma(s, h^{-1}(Y^s))ds,
\end{equation} 
for some $\P$-Brownian motion $W$.
\item Let $W$ be a Brownian motion on $(\Omega, \shf,\P)$.
 $X$ is a solution to the strong martingale
problem with respect to $\shd_L$ and $W$ if and only if
\eqref{Y} holds.

\end{enumerate}
\end{proposition}

\begin{preuve}

 \begin{enumerate}	
\item
We start proving the direct implication.
According to \eqref{lmp} and the notations introduced therein
\begin{equation} \label{Ybis}
M_t^{h} = h(X_t) - h(X_0) - \int_{0}^{t}\mathit{L}h(X_s)ds - \int_{0}^{t}h'(X_s)\Gamma(s, X^s)ds,
\end{equation}
is a $\P$-local martingale
on some probability space $(\Omega, \mathcal{F})$.

In particular, by Corollary \ref{C22}
\begin{equation*}
Y_t = Y_0 + \int_{0}^{t}h'(h^{-1}(Y_s))\Gamma(s, h^{-1}(Y^s))ds + M_t^{h},
\end{equation*}
where $M^h$ is a local martingale, so $Y$ is a semimartingale. We need now evaluate $[M^h]_t = [Y]_t.$
We apply \eqref{lmp} for $f = h^2$ and again by Corollary \ref{C22} we get
\begin{equation}\label{pow2}
Y^{2}_t = Y_0^2 + \int_{0}^{t}\sigma^{2}_0(Y_s)ds + 2\int_{0}^{t}Y_s h'(h^{-1}(Y_s))\Gamma(s, h^{-1}(Y^s))ds + M_t^{h^2},
\end{equation}
where $M^{h^2}$ is a local martingale and
 we recall that $\sigma_0$ was defined in \eqref{Esigma0}.
By integration by parts,
\begin{equation*}
[ Y]_t =  Y^{2}_t -Y_0^2 - 2 \int_{0}^{t}Y_s dY_s = Y^{2}_t -Y_0^2 + M_t -2\int_{0}^{t}Y_s h'(h^{-1}(Y_s))\Gamma(s, h^{-1}(Y^s))ds,
\end{equation*}
where $M_t = -2\int_{0}^{t}Y_{s}dM^{h}_s$.
Therefore
\begin{equation}\label{rpow2}
Y^{2}_t = Y_0^2 -M_t + 2\int_{0}^{t}Y_s h'(h^{-1}(Y_s))\Gamma(s, h^{-1}(Y^s))ds +
[Y]_t.
\end{equation}
Now we can use the uniqueness of the decomposition of a semimartingale $Y^2$
which admits the two expressions
\eqref{pow2} and \eqref{rpow2}. This says $-M = M^{h^2}$ and $\int_{0}^{t}\sigma^{2}_0(Y_s)ds = [Y]_t$.
By \eqref{Ybis}
$$ [M^h]_t = [Y]_t = \int_0^t \sigma_0^2(Y_s) ds.$$
Setting
\begin{equation*}
W_t := \int_{0}^{t}\frac{dM^h_s}{\sigma_{0}(Y_s)}, t \ge 0,
\end{equation*}
we have
$$ [W ]_t \equiv t.$$
Therefore, by L\'{e}vy's characterization of Brownian motion,
$W$ is a standard Brownian motion. Since
$$ M^h = \int_0^\cdot \sigma_{0}(Y_s) dW_s,$$
\eqref{Ybis} shows that
 $Y$ solves \eqref{Y}.

About the converse implication suppose now that $Y = h(X) $ satisfies \eqref{Y},
  for some  $\P$-Brownian motion $W$. We take $f \in \mathcal{D}_{\mathit{L}}$. By
Proposition \ref{P22}
 $g \equiv f\circ h^{-1} \in C^{2}$. Using It\^o's formula,
Proposition \ref{Ph} and
Remark \ref{Remh}  we get
\begin{eqnarray*}
g(Y_t) &=& g(Y_0) + \int_{0}^{t}g'(Y_s)dY_s + \frac{1}{2}\int_{0}^{t}g''(Y_s)d[Y]_s\\ &=&
f(X_0) + \int_{0}^{t}Lf(X_s)ds + \int_{0}^{t}f'(X_s)\sigma(X_s)dW_s + \int_{0}^{t}f'(X_s)\Gamma(s, X^s)ds.
\end{eqnarray*}
Therefore
\begin{equation*}
f(X_t) - f(X_0) - \int_{0}^{t}\mathit{L}f(X_s)ds - \int_{0}^{t}f'(X_s)\Gamma(s, X^s)ds = \int_{0}^{t}f'(X_s)\sigma(X_s)dW_s
\end{equation*}
is a local martingale, which concludes the proof.
\item The converse implication follows in the same way as for item (1).
The proof of the direct implication follows directly by It\^o's formula.
\end{enumerate}
\end{preuve}
\begin{corollary} \label{C32}
Let $(X,\P)$ be a solution to the martingale problem 
related to \eqref{DOperatorL}
with respect to $\shd_L$.
 Then $X$ is a Dirichlet process (with respect to its canonical filtration)
and $[X]_t = \int_0^t \sigma^2(X_s)ds, t \in [0,T]$.
\end{corollary}
\begin{preuve}
	
By Proposition \ref{X-Y}, $X = h^{-1}(Y)$,
where $Y$ is obviously a semimartingale such
that $[Y]_t = \int_0^t \sigma_0^2(Y_s)ds, t \in [0,T].$
Consequently, by Remark \ref{RDir}, $X$ is indeed
a Dirichlet process and
$$[X]_t = \int_0^t ((h^{-1})')^2\sigma^2_{0}(Y_s) ds =
\int_0^t \frac{\sigma_0^2}{(h' \circ h^{-1})^2}(Y_s) ds
= \int_0^t \sigma^2(X_s)ds, t \in [0,T]. $$
\end{preuve}
\begin{remark} \label{R413} If $X$ 
is a solution to the strong martingale problem
with respect to $\shd_L$ and some Brownian motion $W$,
  then $X$ is a Dirichlet process with respect to the canonical filtration
of the Brownian motion.
\end{remark}
An immediate consequence of Proposition \ref{X-Y} is the
following.
\begin{corollary}\label{CorMart}
Suppose that $\Gamma = 0$ and let $(X, \P)$ be a solution of the
martingale problem related to $L$ with respect to $\shd_L$.
Then $Y = h(X)$ is an $\mathfrak{F}$-local martingale
where $\mathfrak{F}$ is the canonical filtration
of $X$ with
quadratic variation
$[Y] = \int_0^\cdot \sigma_0^2(Y_s)ds.$

\end{corollary}

\subsection{Existence}
\label{S43}

$ $

We make here the same conventions as in Section \ref{SPD}.
In the sequel, we introduce the map $\tilde \Gamma: \Lambda \rightarrow \R$
defined by
\begin{equation} \label{GammaTilde}
 \tilde{\Gamma}(s, \eta) = \dfrac{\Gamma(s, \eta)}{\sigma(\eta(s))}, (s, \eta) \in \Lambda.
\end{equation}

At this point, we introduce the following  assumption.
\begin{assumption} \label{H}
  There is $K > 0$ such that
$$
\sup_{s \in [0, T]}|\tilde{\Gamma}(s, h^{-1}\circ \eta^s)| \leq K\left(1 + \sup_{s \in [0, T]}|\eta(s)|\right), \ \forall \eta \in C([0,T]).$$
\end{assumption}
\begin{remark} \label{R417}
	Given a stochastic process $X$, setting $Y = h(X)$, we get \begin{equation}\label{lg}
	\sup_{s \in [0, T]}|\tilde{\Gamma}(s, X^s)| \leq K\left(1 + \sup_{s \in [0, T]}|Y_s|\right).
	\end{equation}
In particular $\int_{0}^{T}\tilde{\Gamma}^2(s, X^s)ds < \infty\ a.s.$
\end{remark}
Proposition \ref{part} below
 is a well-known extension of Novikov's criterion. 
It is an easy consequence of
Corollary 5.14 \cite[Chapter 3]{karatzasShreve}.

\begin{proposition}\label{part} Suppose Assumption \ref{H}.
	Let $W$ be a Brownian motion and $X$ a continuous and adapted process for which there exists a partition $0 = t_0, t_1, ..., t_n = T$ such that for $i \in \{1, ..., n\}$
		$$\E\left[\exp\left(\frac{1}{2}\int_{t_{i - 1}}^{t_i}|\tilde{\Gamma}(s,X^s)|^2 ds\right)\right] < \infty.$$
	Then, the process $$N_t = \exp\left(\int_{0}^{t}\tilde{\Gamma}(s, X^s) dW_s - \frac{1}{2}\int_{0}^{t}|\tilde{\Gamma}(s, X^s)|^2ds\right),$$
	
	is a martingale.
\end{proposition}
We will need a slight adaptation of the Dambis-Dubins-Schwarz theorem to the
case of  a finite interval.
\begin{proposition}\label{DDS}
Let $M$ be a local martingale vanishing at zero
 such that $[M]_t = \int_{0}^{t}A_s ds$, $t \in [0, T].$
 Then, on a possibly enlarged probability space,
 there exists a copy of $M$ (still denoted by the same letter  $M$)
with the same law and a
Brownian motion $\beta$ such that 
$$M_t = \beta_{\int_{0}^{t}A_s ds}, t \in [0,T].$$
\end{proposition}

\begin{preuve}
	
Let us define $$\tilde{M}_t = \left\{\begin{array}{lr}
M_t, & \ t \in [0, T]\\
M_{T} + B_t - B_T, & \ t > T,
\end{array}\right.$$
 where $B$ is a Brownian motion independent of $M$. 
 If the initial probability space is not rich enough, one considers
 an enlarged probability space containing a copy of $M$
 (still denoted by the same letter) with
 the same law and the independent Brownian motion $B$.
Note that $\tilde{M}$ is a local martingale and 
we have, $$[\tilde{M}]_t = \left\{\begin{array}{lr}
[M]_t, & \ t \in [0, T]\\
t - T + [M]_T, & \ t > T.
\end{array}\right.$$
Observe that $\displaystyle\lim_{t \to \infty}[\tilde{M}]_t = \infty$. By the
classical Dambis,
Dubins-Schwarz theorem there exists a 
standard Brownian motion $\beta$ such that a.s.
$\tilde{M}_t = \beta_{\int_{0}^{t}A_s ds}, t \ge 0$.
In particular  $$M_t = \beta_{\int_{0}^{t}A_s ds}, \ 0 \leq t \leq T.$$

\end{preuve}	

The proposition below is an adaptation of a well-known argument
 for Markov diffusions.
\begin{proposition}\label{mart}
Suppose that $\sigma_0$ is bounded.
Let $(X,\P)$ be a solution of
the martingale problem related to \eqref{DOperatorL} with respect to $\shd_L$
with $\Gamma = 0$. Let $M^X$ the martingale component of $X$.
We set
$$ W_t := \int_0^t \frac{1}{\sigma(X_s)}dM^X, t \in [0,T].$$

Then
\begin{equation}\label{expo}\exp\left(\int_{0}^{t}\tilde{\Gamma}(s, X^s) dW_s - \frac{1}{2}\int_{0}^{t}|\tilde{\Gamma}(s, X^s)|^2ds\right), t \in [0,T],
\end{equation}
		is a martingale.
\end{proposition}
\begin{remark}\label{RJust}
We recall that, by Corollary \ref{C32},
 $X$ is an $\mathfrak{F}$-Dirichlet process
($\mathfrak{F}$ be the canonical filtration)
and
$[X] = [M^X] = \int_0^\cdot \sigma^2(X_s)ds$ so that
by L\'evy's characterization theorem, $W$ is an  $\mathfrak{F}$-Brownian motion.

\end{remark}

\begin{prooff} (of Proposition \ref{mart}).
	
	Let $Y = h(X)$. By Proposition \ref{X-Y},
$[Y] = \int_0^\cdot \sigma_0^2(Y_s)ds.$
	Let $k \geq \vert \sigma_0 \vert^2_\infty T$ and $\{t_0 = 0, ..., t_n = T\}$ a grid of $[0, T]$ so that
\begin{equation} \label{Eci}
 c_i:= \frac{3}{2} (t_i - t_{i - 1})K^2 k  < \dfrac{1}{2},
\end{equation}
 ($K$ coming from Assumption \ref{H}) $i = \{1, ..., n\}$. 
By \eqref{lg} we know that
	\begin{equation}\label{lg1}
	\int_{t_{i - 1}}^{t_i}|\tilde{\Gamma}(s, X^s)|^2 ds \leq (t_{i} - t_{i - 1})
K^2\left(1 + \sup_{s \in [0, T]}|Y_s|\right)^2.
	\end{equation}
	We set $M_t = Y_t - Y_0, t \in [0,T]$.
	Note that
\begin{equation} \label{EDDS1}
 \left(1 + \sup_{s \in [0, T]}|Y_s|\right)^2
 \leq 3 \sup_{s \in [0, T]}|M_s|^2 + 3(1 + Y_0^2).
\end{equation}
We recall that $Y_0$ is deterministic.
In view of applying Proposition \ref{part},
taking into account \eqref{lg1} and \eqref{EDDS1} we get
	\begin{eqnarray}\label{lg2}
&&	\E\left(\exp\left(\frac{1}{2}\int_{t_{i - 1}}^{t_i}\tilde{\Gamma}^2
(s, X^s) ds\right)\right)
\leq  \\ 
&& \E \left( \exp\left(\frac{3(t_{i} - t_{i - 1})K^2}{2} 
 \sup_{s \in [0, T]}|M_s|^2\right) \right)
\exp\left(\frac{3}{2} (t_{i} - t_{i - 1})K^2 (1 + Y_0^2)\right). \nonumber
	\end{eqnarray}
	Since $M$ is a local martingale vanishing at zero,
 Proposition \ref{DDS} states that there is a
copy (with the same distribution) of $M$ (still denoted by the same letter)
 on another probability space, a
 Brownian motion $\beta$ such that previous expression gives
\begin{eqnarray} \label{EDDS2}
&& \E \left(\exp\left(\frac{3(t_{i} - t_{i - 1})K^2}{2} 
 \sup_{s \in [0, T]}|\beta_{[M]_s}|^2\right) \right) 
\exp\left(\frac{3}{2} (t_{i} - t_{i - 1})K^2 (1 + Y_0^2)\right) \nonumber\\
& \le& \E \left( \exp\left(\frac{3(t_{i} - t_{i - 1})K^2}{2} 
   \sup_{\tau \in [0, k]}|\beta_{\tau}|^2 
\right) \right) 
\exp\left(\frac{3}{2} (t_{i} - t_{i - 1})K^2 (1 + Y_0^2)\right),
  \end{eqnarray}
 the latter inequality being valid because $[M]_t = \int_0^t \sigma_0^2(Y_s)ds$.
By \eqref{Eci} we get
	\begin{eqnarray} \label{lg2bis}
	\E\left[\displaystyle\exp\left(\frac{1}{2}
\int_{t_{i - 1}}^{t_i}\tilde{\Gamma}^2(s, X^s) ds\right)\right] &\leq&  
\E\left[\exp\left( \frac{c_i}{k} 
\sup_{\tau \in [0,k]}  \vert B_{\tau} \vert ^2)\right) \right]
\exp\left(\frac{c_i}{k}(1+Y_0^2)\right) \nonumber \\
&& \\
&\le& \E\left[ \sup_{\tau \in [0,k]}  \exp\left( \frac{c_i}{ k} 
 \vert B_{\tau} \vert ^2\right) \right]
\exp\left(\frac{c_i}{k}(1+Y_0^2)\right). \nonumber
	\end{eqnarray}
By  Remark \ref{gau} below
$$ \E\left[\exp\left(\frac{c_i}{ k} \vert B_{\tau}\vert ^2\right)\right]  \leq
\E\left[\exp\left(c_i G^2\right)\right] < \infty, $$
where $G$ is a standard Gaussian random variable.
Since $x \mapsto \exp(\frac{c_i}{2k} x)$ is increasing and convex, 
and $ (\vert B_{\tau} \vert^2)_{\tau \ge 0}$ is a non-negative square integrable 
submartingale, 
 then $(\exp(\frac{c_i}{2k} \vert B_{\tau} \vert^2)$ is also a non-negative 
submartingale.
Consequently, by Doob's inequality (with $p=2)$  the expectation on the right-hand side 
of \eqref{lg2bis} is finite.
Finally by Proposition \ref{part}, \eqref{expo} is a martingale.
\end{prooff}
\begin{remark}\label{gau}
	Let $G$ be a standard Gaussian r.v. If $c < \frac{1}{2}$
 then $$\E[\exp(cG^2)] < \infty.$$
\end{remark}
 This opens the way to the following existence result for
our martingale problem.
\begin{theorem} \label{P41}
Suppose the validity of Assumption \ref{A1} and
that one of the two conditions  below are fulfilled.
\begin{enumerate}
\item $\tilde{\Gamma}$ is bounded.
\item $\tilde{\Gamma}$ fulfills Assumption \ref{H} 
and $\sigma_0$ is bounded. 
\end{enumerate}
 Then the martingale problem  related to \eqref{DOperatorL} 
with respect to  $\shd_L$ admits
existence.

\end{theorem}

\begin{preuve}

By Proposition \ref{P21}, we can consider 
 a solution $(X,\P)$  to the above-mentioned martingale problem 
with $\Gamma = 0$. 
By Remark \ref{P32frw1}, there is a Brownian motion $W$ such that
\begin{equation}\label{e2}
f(X_t) - f(X_0) - \int_{0}^{t}Lf(X_s)ds = \int_{0}^{t}(f'\sigma)(X_s)dW_s,
\end{equation}
for every $f \in \shd_L$.
We define  
the process $$V_t := \exp\left(\int_{0}^{t}\tilde{\Gamma}(s, X^s)dW_s - \dfrac{1}{2}\int_{0}^{t}\tilde{\Gamma}^2(s, X^s)ds\right).$$
Under item (1),  $V$ is a martingale by the Novikov condition.
Under item (2), Proposition \ref{mart} says that $V$ is a martingale.
We define
\begin{equation}\label{e1}
\tilde{W}_t:= W_t - \int_{0}^{t}\tilde{\Gamma}(s, X^s)ds.
\end{equation}	
By Girsanov's theorem, 
\eqref{e1} is a Brownian motion under the probability
$\Q$ such that
 $$d\Q:= \displaystyle\exp\left(\int_{0}^{T}\tilde{\Gamma}(s, X^s)dW_s - \dfrac{1}{2}\int_{0}^{T}\tilde{\Gamma}^2(s, X^s)ds\right)d\P.$$
Applying \eqref{e1} in \eqref{e2} we obtain
$$f(X_t) - f(X_0) - \int_{0}^{t}Lf(X_s)ds - \int_{0}^{t}f'(X_s)\Gamma(s, X^s)ds = \int_{0}^{t}(f'\sigma)(X_s)d\tilde{W}_s,$$
for every $f \in \shd_L$.
Since $\int_{0}^{t}(f'\sigma)(X_s)d\tilde{W}_s$ is a local martingale
under $\Q$, $(X,\Q)$ is proved to be a solution to
the martingale problem in the statement.

\end{preuve}

\subsection{Uniqueness in law}
$ $

We use here again the notation $\tilde \Gamma$ introduced in
\eqref{GammaTilde}.
\begin{proposition}\label{uni}
Suppose the validity of
Assumption \ref{A1}.
 Then the martingale problem  related to \eqref{DOperatorL} with respect to $\shd_L$ admits
 uniqueness.
\end{proposition}

\begin{preuve}
	
Let $(X^i, \P^i)$,
$i = 1, 2$, be two solutions of
the martingale problem related to \eqref{DOperatorL} with
respect to $\shd_L$.  Let us fix $i = 1, 2.$
By Corollary \ref{C32}, $X^i$ is a
${\mathfrak F}^{X^i}$-Dirichlet process with respect
to  $\P^i$, 
such that $[X^i] \equiv \int_0 ^\cdot \sigma(X^{i}_s)^2 ds.$
Let
$M^i$ be its respective martingale component.
Since $[M^i] \equiv \int_0 ^\cdot \sigma(X^{i}_s)^2 ds,$
by L\'evy's characterization theorem, the process 
\begin{equation} \label{IntWi}
W^{i}_t = \int_{0}^{t}\dfrac{dM^i_s}{\sigma(X^{i}_s)}, t \in [0,T],
\end{equation}
is an ${\mathfrak F}^{X^i}$-Brownian motion.
In particular, $W^i$ is a Borel functional of $X^i$.

By means of localization (similarly as in Proposition 5.3.10
of \cite{karatzasShreve}, without restriction of generality
we can suppose $\tilde \Gamma$ to be bounded.
We define the process (whose r.v. are also  Borel functionals
of $X^i$)
\begin{equation*}
V^{i}_t
 = \exp\left(- \int_{0}^{t} \tilde{\Gamma}(s, X^{i,s}) dW^{i}_s - \frac{1}{2}\int_{0}^{t}\left(\tilde{\Gamma}(s, X^{i,s})\right)^{2}d s\right),
\end{equation*}
which, by Novikov's condition, is a $\P^i$-martingale. This allows us to
define the probability $Q^i$ defined by $d\Q^i = V_T^i d\P^i$.

By Girsanov's theorem, under
$\Q^{i}$,  $B^i_t := W^{i}_t + \int_{0}^{t}\tilde{\Gamma}(s, X^{i,s})ds$
is a Brownian motion. Therefore,
$(X^i, \Q^i)$ solves
 the martingale problem related  to $L$
($\Gamma = 0$)
with respect to  $\shd_L$.
By uniqueness of the martingale problem with respect to $\shd_L$
and $\Gamma = 0$ (see Proposition \ref{P21}),
 $X^{i}$ (under $\Q^{i}$), $i=1,2$ have the same law.
Hence, for every Borel set $B \in \mathfrak{B}(C[0, T])$, we have
\begin{equation*}
\P^{1}\{X^{1} \in B\} \\=\\ \int_{\Omega}\dfrac{1}{V^1_{T}(X^{1})}1_{\{X^1 \in B\}}d\Q^{1}\\=\\
\int_{\Omega}\dfrac{1}{V^2_{T}(X^{2})}1_{\{X^2 \in B\}}d\Q^{2}\\=\\
\P^{2}\{X^2 \in B\}.
\end{equation*}
Therefore, $X^1$ under  $\P^1$ has the same law as $X^2$ under $\P^2$.
Finally the martingale problem related to \eqref{DOperatorL} with respect
to $\shd_L$ admits uniqueness.
\end{preuve}

\subsection{Results on pathwise uniqueness}
$ $

Before exploring conditions for strong existence and
uniqueness for the martingale problem, we state and
prove Proposition \ref{TStrong}, which constitutes a
crucial preliminary step.

Let $\bar \Gamma: \Lambda \rightarrow \R$ be a generic Borel functional.
Related to it,
we formulate the following.
\begin{assumption} \label{ALip}
	$ $
\begin{enumerate}
\item There exists a function $l: \R_+ \rightarrow \R_+$ such that $\int_{0}^{\epsilon}l^{-2}(u)du = \infty$ for all $\epsilon > 0$ and $$|\sigma_{0}(x) - \sigma_{0}(y)| \leq l(|x - y|).$$
\item
$\sigma_0$ has at most linear  growth.
\item  there exists $K > 0$ such that
  $$ |\bar{\Gamma}(s, \eta^1) - \bar{\Gamma}(s, \eta^2) \vert
 \le K\left(\vert\eta^1(s) - \eta^2(s)  \vert + \int_0^s
   \vert \eta^1(r) - \eta^2(r)\vert dr\right), $$
 for all $(s,\eta^1), (s,\eta^2) \in \Lambda.$
\item $\bar \Gamma_\infty :=\displaystyle\sup_{s \in [0, T]} \vert \bar \Gamma(s, 0)\vert < \infty.$
\end{enumerate}
\end{assumption}

\begin{proposition} \label{TStrong}
Let $y_0 \in \R$.
Suppose the validity
of 
Assumption \ref{ALip}.
Then $E(\sigma_0,0,\bar \Gamma)$, i.e.
 \begin{equation}\label{Ybis1}
Y_t = y_0 + \int_{0}^{t}\sigma_0(Y_s)dW_s + \int_{0}^{t} \bar \Gamma(s, Y^s)ds,
\end{equation}
 admits pathwise uniqueness.
 \end{proposition}
 Before proceeding with the proof of the proposition,
we state a lemma which is an easy consequence
of Problem 5.3.15 of \cite{karatzasShreve}.

\begin{lemma} \label{L1}
Suppose the validity of the assumptions in Proposition \ref{TStrong}.
Let $Y$ be a  solution of
\eqref{Ybis1} and $m \geq 2$ an integer.
Then there exists a constant $C>0$, depending
on the linear growth constant of $\sigma_0$,
$Y_0$, $ T, m$, and the
quantities $K, \bar \Gamma_\infty$ respectively
in Assumptions \ref{ALip} (3)-(4), 
such that
$$ \E\left(\sup_{t \le T} \vert Y_t\vert^m\right) \le C.$$
\end{lemma}
\begin{prooff}  \ (of Proposition \ref{TStrong}).

Now let $Y^1, Y^2$ be two solutions on the same probability space
with respect to the same Brownian motion  $W$ of
 \eqref{Ybis1} such that $Y^1_0 = Y^2_0 = y_0$. By Lemma \ref{L1} 
 we have
\begin{equation} \label{E100}
\E\left(\sup_{t \in [0,T]} \vert Y^i_t \vert^2\right) < \infty,
\end{equation}
for $i = 1,2$. By the assumption on $\sigma_0$, this obviously gives
\begin{equation} \label{E101}
\E\left(\int_0^T\vert \sigma_0(Y^i_t) \vert^2 dt\right) < \infty,
\end{equation}
for $i = 1,2$. We set $\Delta_t = Y^1_t - Y^2_t, t \in [0,T],$ and this gives
\begin{equation}\label{e4}
\Delta_t = \int_{0}^{t}(\bar{\Gamma}(s, Y^{1,s}) - \bar{\Gamma}(s, Y^{2, s}))ds + \int_{0}^{t}(\sigma_{0}(Y^1_s) - \sigma_{0}(Y^2_s))dW_s,  t \in [0,T].
\end{equation}
We recall from the proof of Proposition 2.13 in
\cite[Chapter 5]{karatzasShreve}, the existence of the functions
\begin{equation*}
\Psi_n(x) = \int_{0}^{|x|}\int_{0}^{y}\rho_n(u)dudy,
\end{equation*}
such that for every $ x \in \R$
\begin{equation} \label{EFeaturesPsi}
0 \leq \rho_n(x) \leq \dfrac{2}{nl^2(x)}, \quad |\Psi'_n(x)| \leq 1,
\quad  |\Psi_n(x)| \le |x|, \quad
\displaystyle\lim_{n\to \infty}\Psi_n(x) = |x|.
\end{equation}
By \eqref{e4}, applying It\^o's formula we get
\begin{eqnarray}\label{e5}
\Psi_n(\Delta_t) &=& \int_{0}^{t}\Psi'_n(\Delta_s)[\bar{\Gamma}(s, Y^{1, s}) - \bar{\Gamma}(s, Y^{2, s})]ds + \dfrac{1}{2}\int_{0}^{t}\Psi''_n(\Delta_s)[\sigma_{0}(Y^1_s) - \sigma_{0}(Y^2_s)]^2 ds  \nonumber\\ &+& \int_{0}^{t}\Psi'_n(\Delta_s)[\sigma_{0}(Y^1_s) - \sigma_{0}(Y^2_s)]dW_s.
\end{eqnarray}
Using Assumption \ref{ALip} and  \eqref{EFeaturesPsi} we get
\begin{equation}\label{e6}
\Psi_n(\Delta_t) \leq \int_{0}^{t}K\left(|Y^1_s - Y^2_s| + \int_{0}^{s}|Y^1_r - Y^2_r|dr\right)ds + \dfrac{t}{n} + M_t,
\end{equation}
where $M_t = \int_{0}^{t}\Psi'_n(\Delta_s)[\sigma_{0}(Y^1_s) - \sigma_{0}(Y^2_s)]dW_s$ is a local martingale.
Since $\Psi_n'$ is bounded and by \eqref{E101}, $M$ is a (even square integrable)
martingale. 

We now apply the expectation and the Fubini's theorem in \eqref{e6} to get
\begin{equation}\label{e7}
\E \Psi_n(\Delta_t) \leq K\int_{0}^{t}\E|Y^1_s - Y^2_s|ds + KT\int_{0}^{t}\E|Y^1_s - Y^2_s|ds + \dfrac{t}{n},
\end{equation}
since $\E M_t = 0$. Passing to the limit when $n \rightarrow \infty$,
by Lebesgue's dominated convergence theorem, we get
\begin{equation}
\E|\Delta_t| \leq (K + TK)\int_{0}^{t}\E|\Delta_s|ds,
\end{equation}
so, by the Gronwall's inequality 
 we obtain  $\E|\Delta_t| = 0$. By the continuity of the sample paths of $Y^1, Y^2$ we conclude that $Y^1, Y^2$ are indistinguishable.
\end{prooff}

We come back to the framework of the beginning of Section \ref{S41}.
We suppose again the
validity of Assumption \ref{A1}.
We recall the definition of the harmonic function
$h$ defined by $h(0) = 0, h'(x) = e^{-\Sigma}$, see \eqref{eqh}.
We recall the notations
$\sigma_{0} = (\sigma h') \circ h^{-1}$.
We define
\begin{equation} \label{GammaBar}
\bar{\Gamma}(s, \eta) := h'(h^{-1}(\eta(s)))\Gamma(s, h^{-1}(\eta^s)), \quad
s \in [0,T], \eta \in C([0,T]).
\end{equation}
\begin{corollary} \label{CStrong}
Under Assumptions  \ref{A1}  and \ref{ALip}, the martingale problem
related to \eqref{DOperatorL} with respect to $\shd_L$  admits pathwise uniqueness.
\end{corollary}
\begin{preuve}
\
This follows by Proposition \ref{TStrong} and Proposition \ref{X-Y},
taking into account \eqref{GammaBar}.
\end{preuve}
\begin{theorem}\label{pauni}
Suppose the validity of Assumptions \ref{A1} and \ref{ALip}
related to $\bar \Gamma$ introduced in \eqref{GammaBar}.
Suppose moreover
that one of two hypotheses below are in force.
\begin{enumerate}
\item $\tilde \Gamma$ defined in \eqref{GammaTilde} is bounded.
\item $\sigma_0, \frac{1}{\sigma_0}$ are bounded.

\end{enumerate}
 Then \eqref{Y} admits strong existence and pathwise uniqueness.
\end{theorem}

\begin{preuve}
	
	By  Proposition \ref{TStrong},  pathwise uniqueness holds.
Indeed, by \eqref{GammaBar}, the equation \eqref{Y} is a particular case of \eqref{Ybis1}.

To prove existence we wish to apply Theorem \ref{P41}. 
For this we need to verify that
 either Hypothesis (1) or (2) (in Theorem \ref{P41}) hold.
Hypothesis (1)  in the above statement coincides with Hypothesis (1)
in Theorem \ref{P41}.
Suppose the validity of (2) in the above statement and
we check that
Assumption \ref{H} 
holds true.
By \eqref{GammaTilde} and the definition of $\bar{\Gamma}$ in
\eqref{GammaBar}, we obtain
 \begin{equation}\label{e1bis}
\sigma_{0}(\eta(s))\tilde{\Gamma}(s, (h^{-1}\circ \eta))
 = \bar{\Gamma}(s, \eta).
\end{equation}
By (3) in Assumption \ref{ALip}, $\frac{1}{\sigma_{0}}$ being bounded
there exists a constant $K_1 > 0$ such that
 $$|\tilde{\Gamma}(s, (h^{-1}\circ\eta))| \leq K_1\left(|\eta(s)| + \int_{0}^{s}|\eta(r)|dr\right) + |\bar{\Gamma}(s, 0)|.$$
 By (4) in Assumption \ref{ALip} and the fact that, given a function
 $\gamma \in C(|0,T])$,
 $$\int_0^s \vert \gamma(r)\vert dr \le s \sup_{r\in [0,s]} \vert\gamma (r) \vert,$$
Assumption  \ref{H} follows.

By Theorem \ref{P41},
the martingale problem related to \eqref{DOperatorL} with
respect to $\shd_L$ admits existence and by Proposition \ref{X-Y} (1),
we have that  \eqref{Y} has a (weak) solution.
At this point,  we can apply
Yamada-Watanabe theorem
to guarantee that the solution is actually strong.
We remark that the Yamada-Watanabe theorem (in the path-dependent
case) proof
is the same as the one
 in the Markovian case, which  is for instance stated
in Proposition 3.20 \cite[Chapter 5]{karatzasShreve}.
\end{preuve}	
As a consequence of Proposition \ref{X-Y} and Theorem \ref{pauni}
we obtain the following.

\begin{corollary}\label{pauniMart}
Under the same assumptions as in Theorem \ref{pauni}, 
the  martingale problem related to \eqref{DOperatorL} with respect to
  $\shd_L$, admits
strong existence and pathwise uniqueness.
\end{corollary}

\section{Appendix:   different notions of solutions when
  $b'$ is a function}
\label{Appendix}

$ $

Let us suppose below that
$\sigma, b': \R \rightarrow \R$ are locally bounded Borel functions
and $\Gamma$ as given in \eqref{EGamma}.
As already mentioned, for simplicity we will only consider initial conditions
$x_0$ to be deterministic.

\begin{definition} \label{DOmegaSol}
Let  $(\Omega, \shf, \P)$ be a probabilty space,
$(W_t)_{t \ge 0}$ a Brownian motion and $x_0 \in \R$.
A solution $X$ of \eqref{E1.3} 
(depending on the probability space $(\Omega, \shf, \P)$,
the Brownian motion $W$ 
and an initial condition $x_0$)
is a progressively measurable process, with respect to $\mathfrak F^W$,
fulfilling \eqref{E1.3}. 
That equation   \eqref{E1.3} will be denoted by $E(\sigma,b', \Gamma)$
 (without specification of the initial condition).
\end{definition}
\begin{definition} \label{D112} ({\bf Strong existence}).

We will say that equation $E(\sigma,b', \Gamma)$ admits {\bf strong existence}
if the following holds.
Given  any probability space
$(\Omega, \shf, \P)$, 
a 
Brownian motion 
$(W_t)_{t \ge 0}$ and 
 $x_0 \in \R$, there exists a process $(X_t)_{t \ge 0}$ which is
solution to $ E(\sigma,b',\Gamma)$  with $X_0 = x_0 $ a.s.
\end{definition}
\begin{definition} \label{D113} ({\bf Pathwise uniqueness}).
We will say that equation $E(\sigma,b', \Gamma)$  admits pathwise uniqueness
if the following property is fulfilled.

Let $(\Omega, \shf, \P)$ be a probability space
carrying a Brownian motion
$(W_t)_{t \ge 0}$.
If two processes $X, \tilde X$ are two solutions
to $E(\sigma,b', \Gamma)$ 
 such that $X_0 = \tilde X_0$ a.s., then $X$ and $\tilde X$
are indistinguishable.
\end{definition}
\begin{definition} \label{D114} ({\bf Existence in law} or {\bf weak existence}).
Let $\nu$ be a probability law on $\R$.
We will say that $E(\sigma,b', \Gamma;\nu)$ admits weak existence if there exists
a probability space $(\Omega, \shf, \P)$ carrying
a Brownian motion $(W_t)_{t \ge 0}$ and a process $(X_t)_{t \ge 0}$ 
such that $(X,\P)$ is a (weak) solution of
$E(\sigma,b',\Gamma;\nu)$, see Definition \ref{D112Bis}. 

We say that $ E(\sigma,b',\Gamma)$ admits weak existence if $E(\sigma,b', \Gamma; \nu)$  admits weak existence
for every $\nu$. 
\end{definition}
\begin{definition} \label{D115} ({\bf Uniqueness in law}).
 Let $\nu$ be a probability law on $\R$. We say that $E(\sigma,b', \Gamma; \nu)$
has a {\bf unique solution in law} if the following holds.
Suppose we have a probability space
$(\Omega, \shf, \P)$ (respectively $(\tilde \Omega, \tilde \shf, \tilde \P)$)
carrying a Brownian motion
$(W_t)_{t \ge 0}$ (respectively $(\tilde W_t)_{t \ge 0}$).
We suppose that
a  process $(X_t)_{t \ge 0}$  (resp. a process $ (\tilde X_t)_{t \ge 0}$)
is a solution of
$E(\sigma,b', \Gamma)$
 such that both   $X_0$  and $\tilde X_0$ are distributed according to  $\nu$.
 Uniqueness in law means that $X$ and $\tilde X$ must have the same law as
 random elements taking  values in
 $ C([0,T])$ or  $C(\R_+).$

We say that $E(\sigma,b',\Gamma)$  has a unique solution in law if $E(\sigma,b', \Gamma; 
\nu)$ has a unique
solution in law for every $\nu$.
  \end{definition}

  \vfill \eject

 {\bf ACKNOWLEDGEMENTS.}
The research related to this paper
was financially supported by the Regional Program MATH-AmSud 2018
grant 88887.197425/2018-00. A.O.
acknowledges the financial support of CNPq Bolsa de
Produtividade de Pesquisa grant 303443/2018-9.

\bibliographystyle{plain}
\bibliography{../../../../BIBLIO_FILE/biblio-PhD-Alan}
\end{document}